\documentclass{article}

\title{Optimal strategy in the game Risk or Safety}
\author{Rüdiger Jehn\footnote{rjehn@yahoo.com} }
\date{March 2026}

\usepackage[utf8]{inputenc}
\usepackage{multirow}
\usepackage{tikz}
\usetikzlibrary{calc, patterns}
\usepackage{amsmath}
\usepackage{lmodern}

\usepackage{color, colortbl}
\usepackage[colorlinks,citecolor=blue,urlcolor=blue,bookmarks=false,hypertexnames=true]{hyperref} %\definecolor{name}{system}{definition}
\definecolor{Gray}{gray}{0.9}

\begin{document}

\maketitle

{\bf Abstract}\\

In the paper it is proven that the two-players turn-based stochastic game "Risk or Safety" has a unique solution. Both players need to play the same strategy if they want to maximize their winning chances. An analytical method based on the transition equations depending on the players decisions at each game situation is presented. However, the method requires the solution of many large linear equation systems, which makes the problem practically unsolvable for more than 6 coins. Hence, an iterative method is used to determine the optimal strategy for up to 20 coins. A look-up table gives the number of coins a player should toss at each possible game situation. 

\section{Introduction}

The game "Risk or Safety" consists of two competing players. A player whose turn it is tosses a coin. If it is heads they earn a point, and they can choose to go for safety, put the point aside in a save box, and give the turn to their opponent. Or they take the risk to continue and toss the coin again. If it is heads again, the number of "open" points increases by one, and they can choose again to continue or stop, converting the "open" points to "saved" points. Whenever it is tails, all the open points are lost, and it is the turn of the other player. Who collects $n$ points first wins.\\

The obvious question is, when is it recommended to take the risk to continue and when to stop and secure the points in hand. Two ways to determine the optimal strategy were applied. The analytic method works only until $n=6$ because the number of possible strategies that need to be analyzed grows rapidly. Therefore for $n=7$ to 20 the results were obtained by simulating the game iteratively, selecting the solution with the highest winning probability after $k$ turns, and continuing this process until no further improvements are made.\\

One variant of the game is, that the player whose turn it is chooses a number $x \le n$ coins, tosses them and if all are heads, they earn $x$ points. If at least one coin is tails, they earn zero points. The player take turns. Logically this is exactly the same game, although psychologically tossing the coins individually and taking a decision after each throw may increase the suspense.  

\section{Uniqueness of solution}

Like the game "Super Six" \cite{jehn2021b}, this game is a two-player turn-based stochastic game (TBSG). Cui and Yang \cite{cui} state that as a zero-sum game, TBSG is known to have a Nash equilibrium strategy, which means there exist a strategy pair that neither agent can benefit from unilaterally deviating. As reference for this statement, they cite the original work by Shapely\cite{shapely}.\\

In Chapter 3 we want to exploit the fact, that the strategy pair is symmetric, i.\,e.\ that both players apply the same optimal strategy. In a paper by Glenn et al.~\cite{glenn}, the two-player game "Can't Stop" is analyzed. The principle of the game is the same as in our game where after each successful throw the player needs to decide to continue or to stop, but the number of possible game situations is orders of magnitude larger. Theorem 2 in this paper states that the game graph of this game has a unique solution, where for all positions a corresponding winning probability in the interval [0, 1] exists.\\

In the proof of Theorem 2 the game graph is split into strongly connected components, and the components are considered in bottom-up order. We demonstrate this approach for our game with 3 coins. Bottom-up means we have first solved the problem for 1 coin (this is trivial: the first player wins with a probability of 2/3) and in the next step for 2 coins (we omit this step).\\

A game position is given by a triple $(a_0, a_1, b)$, where $a_0$ and $a_1$ are the open and saved points of player~1, respectively, and $b$ is the number of saved points of player~2. We start with the situation (0, 2, 0), i.~e.\ player 1 has already two points saved. Let $x_2$ be the chance of player~1 winning when player~1 starts and $y_2$ the chance of player~1 winning when player~2 starts. If player~1 tosses heads they win, if they toss tails, the turn switches to player~2, hence
\begin{equation}
   x_2 = \frac{1}{2} + \frac{1}{2} y_2.
\end{equation}
On the other hand, if it is player's 2 turn, there are three possibilities: Player 2 stops after
\begin{enumerate}
    \item throw, then $y_2 = \frac{1}{2} x_2 + \frac{1}{2} P(0, 2, 1) = \frac{1}{2} x_2 + \frac{2}{5}$ 
    \item throw, then $y_2 = \frac{3}{4} x_2 + \frac{1}{4} P(0, 2, 2) = \frac{3}{4} x_2 + \frac{1}{6}$ 
    \item throw, then $y_2 = \frac{7}{8} x_2$ 
\end{enumerate}
where $P(a_0, a_1, b)$ denotes the probability for player~1 to win at position $(a_0, a_1, b)$. In previous steps we have calculated $P(0, 2, 1) = \frac{4}{5}$ and $P(0, 2, 2) = \frac{2}{3}$. Since player 2 tries to minimize the winning chances of player 1
\begin{equation}
    y_2 = \min \left ( \frac{1}{2} x_2 + \frac{2}{5}, \frac{3}{4} x_2 + \frac{1}{6}, \frac{7}{8} x_2 \right ).
\end{equation}

Equations 1 and 2 can be iteratively solved (starting for instance with $x_2=0$) and yields the results $x_2 = P(0, 2, 0) = \frac{8}{9}$ and $y_2 = \frac{7}{9}$.\\

In the next step we determine $P(0, 1, 0)$.  Let $x_1$ be the chance of winning when player~1 starts and $y_1$ the chance of winning when player~2 starts. Player~1 has two possibilities, they can stop after
\begin{enumerate}
    \item throw, then $x_1 = \frac{1}{2} y_1 + \frac{1}{2} y_2 = \frac{1}{2} y_1 + \frac{7}{18}$ 
    \item throw, then $x_1 = \frac{3}{4} y_1 + \frac{1}{4}$ 
\end{enumerate}

Since the player wants to maximize their winning chances
\begin{equation}
   x_1 = \max \left ( \frac{1}{2} y_1 + \frac{7}{18}, \frac{3}{4} y_1 + \frac{1}{4} \right ).
\end{equation}
On the other hand, if it is player's 2 turn, there are three possibilities: Player 2 stops after
\begin{enumerate}
    \item throw, then $y_1 = \frac{1}{2} x_1 + \frac{1}{2} P(0, 1, 1) = \frac{1}{2} x_1 + \frac{2}{7}$ 
    \item throw, then $y_1 = \frac{3}{4} x_1 + \frac{1}{4} P(0, 1, 2) = \frac{3}{4} x_1 + \frac{1}{10}$ 
    \item throw, then $y_1 = \frac{7}{8} x_1$ 
\end{enumerate}
In previous steps we have calculated $P(0, 1, 1) = \frac{4}{7}$ and $P(0, 1, 2) = \frac{2}{5}$. Since player 2 tries to minimize the winning chances of player 1
\begin{equation}
    y_1 = \min \left ( \frac{1}{2} x_1 + \frac{2}{7}, \frac{3}{4} x_1 + \frac{1}{10}, \frac{7}{8} x_1 \right ).
\end{equation}

Equations 3 and 4 can be iteratively solved and yields the results $x_1 = P(0, 1, 0) = \frac{8}{11}$ and $y_1 = \frac{7}{11}$.\\

The last step for the game with three coins is the determination of $P(0, 0, 0)$.  Let $x_0$ be the chance of winning when player~1 starts and $y_0$ the chance of winning when player~2 starts. Player~1 has three possibilities, they can stop after
\begin{enumerate}
    \item throw, then $x_0 = \frac{1}{2} y_0 + \frac{1}{2} y_1 = \frac{1}{2} y_0 + \frac{7}{22}$ 
    \item throw, then $x_0 = \frac{3}{4} y_0 + \frac{1}{4} y_2 = \frac{3}{4} y_0 + \frac{7}{36}$ 
    \item throw, then $x_0 = \frac{7}{8} y_0 + \frac{1}{8}$ 
\end{enumerate}

Since the player wants to maximize their winning chances
\begin{equation}
   x_0 = \max \left ( \frac{1}{2} y_0 + \frac{7}{22}, \frac{3}{4} y_0 + \frac{7}{36},  \frac{7}{8} y_0 + \frac{1}{8} \right ).
\end{equation}
On the other hand, if it is player's 2 turn, there are three possibilities: Player 2 stops after
\begin{enumerate}
    \item throw, then $y_0 = \frac{1}{2} x_0 + \frac{1}{2} P(0, 0, 1) = \frac{1}{2} x_0 + \frac{2}{11}$ 
    \item throw, then $y_0 = \frac{3}{4} x_0 + \frac{1}{4} P(0, 0, 2) = \frac{3}{4} x_0 + \frac{1}{18}$ 
    \item throw, then $y_0 = \frac{7}{8} x_0$ 
\end{enumerate}
In previous steps we have calculated $P(0, 0, 1) = \frac{4}{11}$ and $P(0, 0, 2) = \frac{2}{9}$. Since player 2 tries to minimize the winning chances of player 1
\begin{equation}
    y_0 = \min \left ( \frac{1}{2} x_0 + \frac{2}{11}, \frac{3}{4} x_0 + \frac{1}{18}, \frac{7}{8} x_0 \right ).
\end{equation}

Equations 5 and 6 can be iteratively solved and yields the results $x_0 = P(0, 0, 0) = \frac{6}{11}$ and $y_0 = \frac{5}{11}$.\\

What we have demonstrated here for the game with 3 coins can be applied in the next step for 4 coins and so on. For any given game situation the optimal strategy and the corresponding winning probabilities can be determined. They are unique and they are valid for both players, hence both players have to apply the same strategy if they want to maximize their winning chances.

%Let $(x^*, y^*)$ denote the Nash equilibrium strategy pair, with $x^*$ being the optimal strategy for player~A and $y^*$ the optimal strategy for player~B. All game situations which payer~A can reach, also player~B can reach, hence we have a symmetric game and $(y^*, x^*)$ is the Nash equilibrium strategy pair if the two player have swapped the positions (e.g.~if player~A starts with tossing tails). Therefore we can conclude that $x^* = y^*$.\\

\section{Analytic method to determine optimal strategy}

We are no longer differentiating between player 1 and 2, but only between the player whose turn it is to toss the coin and the opponent. As above, each game position is described by the triple $(a_0, a_1, b)$ where $a_0$ are the open points and $a_1$ are the saved points of the player whose turn it is to toss the coin and $b$ are the saved points of the opponent. 

\subsection{Two points to win}
We start with the analysis when two points are needed to win the game. The position (0, 1, 1) is trivial: the probability to win for the player tossing the coin is $\frac{1}{2}+ \frac{1}{8} + \frac{1}{32} + \ldots = \frac{2}{3}$.  For the other 5 possible positions the probabilities depend on the chosen strategies of the players. This diagram depicts all possible game ramifications.\\

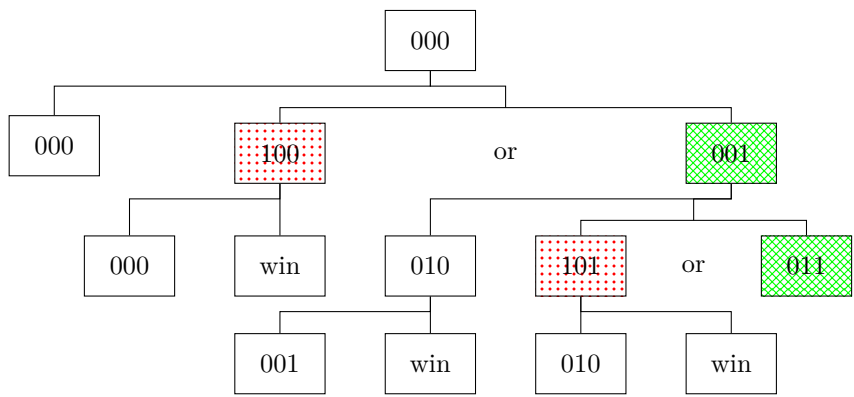
\begin{tikzpicture}[
  every node/.style={draw, minimum width=1.2cm, minimum height=0.8cm, align=center}]

% nodes
\node (root) at (3,0) {000};
\node (left) at (-2,-1.4) {000};
\node[pattern=dots, pattern color=red] (right1) at (1,-1.5) {100};
\node[pattern=crosshatch, pattern color=green] (right2) at (7,-1.5) {001};
\node (21) at (-1,-3) {000};
\node (22) at (1,-3) {win};
\node (23) at (3,-3) {010};
\node[pattern=dots, pattern color=red] (24) at (5,-3) {101};
\node[pattern=crosshatch, pattern color=green] (25) at (8,-3) {011};
\node (31) at (1,-4.3) {001};
\node (32) at (3,-4.3) {win};
\node (33) at (5,-4.3) {010};
\node (34) at (7,-4.3) {win};

% connector lines
\draw (root.south) -- ++(0,-0.2) -| (left.north);
\draw (root.south) -- ++(0,-0.2) -| ($ (right1.north)!0.5!(right2.north) + (0,0.2) $);
\draw (right1.north) -- ++(0,0.2) -| (right2.north);
\draw (right1.south) -- ++(0,-0.2) -| (21.north);
\draw (right1.south) -- ++(0,-0.2) -| (22.north);
\draw (right2.south) -- ++(0,-0.2) -| (23.north);
\draw (right2.south) -- ++(0,-0.2) -| ($ (24.north)!0.5!(25.north) + (0,0.2) $);
\draw (24.north) -- ++(0,0.2) -| (25.north);
\draw (23.south) -- ++(0,-0.2) -| (31.north);
\draw (23.south) -- ++(0,-0.2) -| (32.north);
\draw (24.south) -- ++(0,-0.2) -| (33.north);
\draw (24.south) -- ++(0,-0.2) -| (34.north);
% "or" labels
\node[draw=none] (1x) at ($(right1)!0.5!(right2)$) {or};
\node[draw=none] (2x) at ($(24)!0.5!(25)$) {or};
\end{tikzpicture}\\

If tails is tossed, i.e.\ no point, the left branches apply. If at position (0, 0, 0) heads is tossed, the right branch applies and if the player chooses "risk", we go to position (1, 0, 0) which is colored with red dots. If, however, the player chooses "safety", they convert this one point and since the opponent takes over, the new game situation is (0, 0, 1) which is colored green.\\

Let $ P_{a_0a_1b} = P(a_0, a_1, b)$ denote the probability to win at position $(a_0, a_1, b)$, then the following equations hold:
$$P_{000}=\frac{1}{2}(1-P_{000}) + \frac{1}{2}P_{100} \text{ or}$$
$$P_{000}=\frac{1}{2}(1-P_{000}) + \frac{1}{2}(1-P_{001})$$
The first equation holds if the player takes the risk and continues to toss the coin and the second if they stop and handover the game to the opponent.\\

Similarly we get
$$P_{001}=\frac{1}{2}(1-P_{010}) + \frac{1}{2}P_{101} \text{ or}$$
$$P_{001}=\frac{1}{2}(1-P_{010}) + \frac{1}{2}(1-P_{011})$$
with $P_{011} = \frac{2}{3}$.\\

For the other three positions, there is just one equation:
$$P_{010}= \frac{1}{2}(1-P_{001}) + \frac{1}{2}$$
$$P_{100}= \frac{1}{2}(1-P_{000}) + \frac{1}{2}$$
$$P_{101}= \frac{1}{2}(1-P_{010}) + \frac{1}{2}$$

There are three different strategies: 
\begin{enumerate}
    \item Continue at (1, 0, 0),
    \item Stop at (1, 0, 0) and continue at (1, 0, 1) and
    \item Stop at (1, 0, 0) and stop at (1, 0, 1).
\end{enumerate}  
Hence we get three different systems of 5 linear equations depending on the three possible strategies. How do we determine effectively the best strategy?\\

If the players stop at (1, 0, 0), the winning probability at position (0, 0, 1) is $\frac{2}{5}$ if the players continue at (1, 0, 1) and $\frac{2}{9}$ if the players stop at (1, 0, 1). Therefore strategy 2 dominates strategy 3. In the next step we compare the winning probability at position (0, 0, 0). For strategy 2 it is $\frac{8}{15}$ which is less than for strategy 1 where a value of $\frac{4}{7}$ is obtained. Hence strategy 1 is optimal.\\

\subsection{Three points to win}

There are 12 new game positions for which the equations need to be established:
$$P_{000}=\frac{1}{2}(1-P_{000}) + \frac{1}{2}P_{100} \text{ or } P_{000}=\frac{1}{2}(1-P_{000}) + \frac{1}{2}(1-P_{001})$$
$$P_{001}=\frac{1}{2}(1-P_{010}) + \frac{1}{2}P_{101} \text{ or }P_{001}=\frac{1}{2}(1-P_{010}) + \frac{1}{2}(1-P_{011})$$
$$P_{010}=\frac{1}{2}(1-P_{001}) + \frac{1}{2}P_{110} \text{ or }P_{010}=\frac{1}{2}(1-P_{001}) + \frac{1}{2}(1-P_{002})$$
$$P_{100}=\frac{1}{2}(1-P_{000}) + \frac{1}{2}P_{200} \text{ or }P_{100}=\frac{1}{2}(1-P_{000}) + \frac{1}{2}(1-P_{001})$$
$$P_{002}=\frac{1}{2}(1-P_{020}) + \frac{1}{2}P_{102} \text{ or }P_{002}=\frac{1}{2}(1-P_{020}) + \frac{1}{2}(1-P_{021})$$
$$P_{020}=\frac{1}{2}(1-P_{002}) + \frac{1}{2}$$
$$P_{101}=\frac{1}{2}(1-P_{010}) + \frac{1}{2}P_{201} \text{ or }P_{101}=\frac{1}{2}(1-P_{010}) + \frac{1}{2}(1-P_{011})$$
$$P_{110}=\frac{1}{2}(1-P_{001}) + \frac{1}{2}$$
$$P_{200}=\frac{1}{2}(1-P_{000}) + \frac{1}{2}$$
$$P_{102}=\frac{1}{2}(1-P_{020}) + \frac{1}{2}P_{202} \text{ or }P_{102}=\frac{1}{2}(1-P_{020}) + \frac{1}{2}(1-P_{021})$$
$$P_{201}=\frac{1}{2}(1-P_{010}) + \frac{1}{2}$$
$$P_{202}=\frac{1}{2}(1-P_{020}) + \frac{1}{2}$$\\

The left branch of the tree where we continue at position (1, 0, 0) shows 5 possible strategies:\\ 

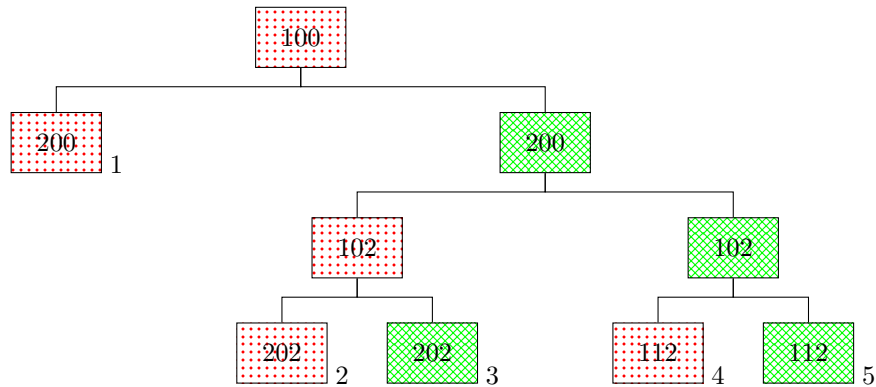
\begin{tikzpicture}[
  every node/.style={draw, minimum width=1.2cm, minimum height=0.8cm, align=center}]

% nodes
\node[pattern=dots, pattern color=red] (100) at (2.25,0) {100};
\node[pattern=dots, pattern color=red] (200) at (-1,-1.4) {200};
\node[pattern=crosshatch, pattern color=green] (002) at (5.5,-1.4) {200};
\node[pattern=dots, pattern color=red] (102) at (3,-2.8) {102};
\node[pattern=crosshatch, pattern color=green] (012) at (8,-2.8) {102}; 
\node[pattern=dots, pattern color=red] (202) at (2,-4.2) {202};
\node[pattern=crosshatch, pattern color=green] (022) at (4,-4.2) {202}; 
\node[pattern=dots, pattern color=red] (112) at (7,-4.2) {112};
\node[pattern=crosshatch, pattern color=green] (022b) at (9,-4.2) {112}; 

% connector lines
\draw (100.south) -- ++(0,-0.25) -| (200.north);
\draw (100.south) -- ++(0,-0.25) -| (002.north);
\draw (002.south) -- ++(0,-0.25) -| (012.north);
\draw (002.south) -- ++(0,-0.25) -| (102.north);
\draw (102.south) -- ++(0,-0.25) -| (202.north);
\draw (102.south) -- ++(0,-0.25) -| (022.north);
\draw (012.south) -- ++(0,-0.25) -| (112.north);
\draw (012.south) -- ++(0,-0.25) -| (022b.north);

\node[draw=none] (1x) at ($(200)+(0.8,-0.3)$) {1};
\node[draw=none] (2x) at ($(202)+(0.8,-0.3)$) {2};
\node[draw=none] (3x) at ($(022)+(0.8,-0.3)$) {3};
\node[draw=none] (4x) at ($(112)+(0.8,-0.3)$) {4};
\node[draw=none] (5x) at ($(022b)+(0.8,-0.3)$) {5};

\end{tikzpicture}\\

Again, the red dotted boxes denote the risk strategy and the green boxes the safety strategies. The numbers next to the boxes denote the strategy number. The comparison of strategy 4 and 5 is not necessary anymore because we know from $n=2$ of the previous section that strategy 4 is superior. In order to determine the winner between strategy 2 and 3 we need to compare the winning probability at (1, 0, 2) which is $\frac{1}{3}$ if the player continues after a successful throw (strategy 2) and only $\frac{1}{5}$ if the player stops (strategy 3). Hence strategy 2 dominates strategy 3.\\

In the next step we compare the winning probability at (0, 0, 2) which is $\frac{2}{9}$ if the player continues after a successful throw (strategy 2) and only $\frac{2}{15}$ if the player stops (strategy 4). Hence strategy 2 dominates strategy 4.\\

In the last step we compare the winning probability at (1, 0, 0) which is $\frac{3}{5}$ if the player continues after a successful throw (strategy 1) and $\frac{13}{21}$ if the player stops (strategy 2) which is slightly larger. Hence strategy 2 dominates strategy 1.\\

The same exercise needs to be done on the right branch of the tree, where in total 26~strategies exist. If we stop at position (1, 0, 0), the players will need to take the next decision either in position (1, 0, 1) or in (1, 1, 0). Hence we have 4~cases to study. This is illustrated in the second line of this diagram:\\

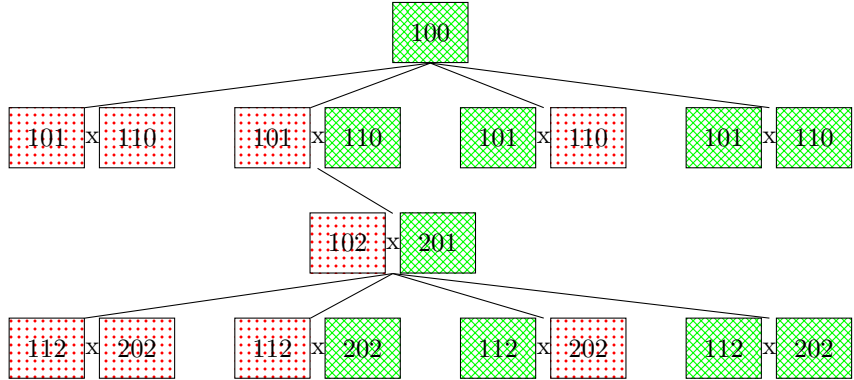
\begin{tikzpicture}[
  every node/.style={draw, minimum width=1.0cm, minimum height=0.8cm, align=center}]

% nodes
\node[pattern=crosshatch, pattern color=green] (1) at (3.1,0) {100};

\node[pattern=dots, pattern color=red] (21) at (-2,-1.4) {101};
\node[pattern=dots, pattern color=red] (22) at (-0.8,-1.4) {110};
\node[draw=none] (2x) at ($(21.east) + (.1,0)$) {x};
\node[pattern=dots, pattern color=red] (23) at (1,-1.4) {101};
\node[pattern=crosshatch, pattern color=green] (24) at (2.2,-1.4) {110};
\node[draw=none] (3x) at ($(23.east) + (.1,0)$) {x};
\node[pattern=crosshatch, pattern color=green] (25) at (4,-1.4) {101};
\node[pattern=dots, pattern color=red] (26) at (5.2,-1.4) {110};
\node[draw=none] (4x) at ($(25.east) + (.1,0)$) {x};
\node[pattern=crosshatch, pattern color=green] (27) at (7,-1.4) {101};
\node[pattern=crosshatch, pattern color=green] (28) at (8.2,-1.4) {110};
\node[draw=none] (5x) at ($(27.east) + (.1,0)$) {x};

\node[pattern=dots, pattern color=red] (31) at (2,-2.8) {102};
\node[pattern=crosshatch, pattern color=green] (32) at (3.2,-2.8) {201};
\node[draw=none] (5x) at ($(31.east) + (.1,0)$) {x};

\node[pattern=dots, pattern color=red] (41) at (-2,-4.2) {112};
\node[pattern=dots, pattern color=red] (42) at (-0.8,-4.2) {202};
\node[draw=none] (2x) at ($(41.east) + (.1,0)$) {x};
\node[pattern=dots, pattern color=red] (43) at (1,-4.2) {112};
\node[pattern=crosshatch, pattern color=green] (44) at (2.2,-4.2) {202};
\node[draw=none] (3x) at ($(43.east) + (.1,0)$) {x};
\node[pattern=crosshatch, pattern color=green] (45) at (4,-4.2) {112};
\node[pattern=dots, pattern color=red] (46) at (5.2,-4.2) {202};
\node[draw=none] (4x) at ($(45.east) + (.1,0)$) {x};
\node[pattern=crosshatch, pattern color=green] (47) at (7,-4.2) {112};
\node[pattern=crosshatch, pattern color=green] (48) at (8.2,-4.2) {202};
\node[draw=none] (5x) at ($(47.east) + (.1,0)$) {x};

% connector lines
\draw (1.south) -- ($(21.east) + (0,0.4)$);
\draw (1.south) -- ($(23.east) + (0,0.4)$);
\draw (1.south) -- ($(25.east) + (0.1,0.4)$);
\draw (1.south) -- ($(27.east) + (0.1,0.4)$);

\draw ($(23.south) + (0.6,0)$) -- ++(1,-0.6);

\draw ($(31.south) + (0.6,0)$) -- ($(41.east) + (0,0.4)$);
\draw ($(31.south) + (0.6,0)$) -- ($(43.east) + (0,0.4)$);
\draw ($(31.south) + (0.6,0)$) -- ($(45.east) + (0.1,0.4)$);
\draw ($(31.south) + (0.6,0)$) -- ($(47.east) + (0.1,0.4)$);

%\draw (100.south) -- ++(0,-0.25) -| (002.north);

\end{tikzpicture}\\

We pick only the case where the players decide at position (1, 0, 1) to continue (=red) and at position (1, 1, 0) to stop (=green). From this branch the players can reach the positions (2, 0, 1) and (1, 0, 2) where new decisions are required. Again there are 4 possible combinations of stop and continue and we pick for the illustration "continue" at (1, 0, 2) and "stop" at (2, 0, 1). Again 4 options open up: stop or continue at the positions (1, 1, 2) and (2, 0, 2). It should be noted that in the implementation of this approach, the option to stop at (1, 1, 2) is not investigated anymore, because we know from $n=2$ that this is a bad strategy. Hence, for this branch, only the two cases stop or continue at position (2, 0, 2) need to be compared. As it turns out, "continue" is the optimal strategy at this situation.\\

If we reuse the results from the case $n=2$ out of the 31 possible strategies, 16 remain, i.e.\ 16 times a $12 \times 12$ matrix is inverted and the results are compared. The optimal strategy is to continue in all situations except at position (1, 0, 0) where it is recommended to stop. This yields a winning probability of $\frac{6}{11}$ ($54.54~\%$) for the player who starts.

\subsection{Four and more points to win}
In case of $n=4$, there are 22 new game situations and out of 3011 total strategies 125 need to be compared. The optimal strategy is to stop at (1, 2, 0), (2, 0, 0) and (2, 0, 2) and else to continue. This yields a winning probability of $\frac{2236}{4165}$ ($53.69~\%$) for the player who starts.\\

In case of $n=5$, there are 35 new game situations and out of 5,755,251 total strategies 1296 need to be compared. The optimal strategy yields a winning probability of $\frac{1026}{1925}$ ($53.30~\%$) for the player who starts.\\

In case of $n=6$, there are 51 new game situations and out of 357 billion possible strategies 16807 need to be compared. The optimal strategy yields a winning probability of $\frac{275848876}{521145625}$ ($52.93~\%$) for the player who starts.\\

In general there are $\frac{n}{2}(3n - 1)$ new game situations, which is the dimension of the matrices that need to be inverted for each possible strategy. Therefore, this analytical approach very soon reaches a numerical limit.

\section{Iterative method to determine optimal strategy}

%let P(a,b,c) denote the probability to win at position (a,b,c). And we can define P(a,b,c) = 1 if a+b=n.
For any $b$ and $c$, we have $P(0,b,c) = \frac{1}{2} (1 - P(0,c,b)) + \frac{1}{2} P(1,b,c)$ and $P(0,c,b) = \frac{1}{2} (1 - P(0,b,c)) + \frac{1}{2} P(1,c,b)$.
Inserting the second equation into the first yields
\begin{equation}
    P(0,b,c) = \frac{1}{3} + \frac{2}{3} P(1,b,c) - \frac{1}{3} P(1,c,b)
    \label{eq:1}
\end{equation}

Applying Eq.~\ref{eq:1} at position $(a,b,c)$ for $a \ge 1$ we get:

$$P(a,b,c) = \frac{1}{2} (1 - P(0,c,b)) + \frac{1}{2} P(a+1,b,c)$$ 
$$= \frac{1}{3} - \frac{1}{3} P(1,c,b) + \frac{1}{6} P(1,b,c) + \frac{1}{2} P(a+1,b,c)$$
when the current player chooses to continue and 
$$P(a,b,c) = 1 - P(0,c,a+b) = \frac{2}{3} - \frac{2}{3} P(1,c,a+b) + \frac{1}{3} P(1,a+b,c)$$
when the current player chooses to stop.\\

We only need to iterate over the positions $(a,b,c)$ with $0 < a < n$,  $0 < a+b < n$ and $0 \le c < n$. There are $\frac{(n-1)n^2}{2}$ such positions.\\

For every position $(a,b,c)$, let $P_{min}(a,b,c)$ be the currently known minimum winning probability and $P_{max}(a,b,c)$ the maximum winning probability obtained in the iteration. And let $S(a,b,c)$ be the current player's strategy with $S(a,b,c) \in \{ \text{unknown, continue, stop} \}$.\\

The first step is the initialisation of all varaibles:
For every position $(a,b,c)$, $P_{min}(a,b,c) := 0$, $P_{max}(a,b,c) := 1$ and $S(a,b,c)$ := unknown.\\

During the iteration steps the three possible cases for $S$ are differentiated:
\begin{enumerate}
  \item $S(a,b,c) =$ unknown. \\[4pt]
    If 
   $\frac{1}{3} - \frac{1}{3} P_{\max}(1,c,b) + \frac{1}{6} P_{\min}(1,b,c) + \frac{1}{2} P_{\min}(a+1,b,c)
    > 
    \frac{2}{3} - \frac{2}{3} P_{\min}(1,c,a+b) + \frac{1}{3} P_{\max}(1,a+b,c),$
    then $S(a,b,c) =$ continue. 

    \medskip
    If 
    $\frac{2}{3} - \frac{2}{3} P_{\max}(1,c,a+b) + \frac{1}{3} P_{\min}(1,a+b,c)
    >
    \frac{1}{3} - \frac{1}{3} P_{\min}(1,c,b) + \frac{1}{6} P_{\max}(1,b,c) + \frac{1}{2} P_{\max}(a+1,b,c),$
    then $S(a,b,c) =$ stop.

  \begin{align}
  P_{\min}(a,b,c) := \min\Big(
      &\tfrac{1}{3} - \tfrac{1}{3} P_{\max}(1,c,b)
      + \tfrac{1}{6} P_{\min}(1,b,c)
      + \tfrac{1}{2} P_{\min}(a+1,b,c), \nonumber\\
      &\tfrac{2}{3} - \tfrac{2}{3} P_{\max}(1,c,a+b)
      + \tfrac{1}{3} P_{\min}(1,a+b,c)
  \Big). \nonumber
  \end{align}

  \begin{align}
  P_{max}(a,b,c) := \max\Big( 
    &\tfrac{1}{3} - \tfrac{1}{3} P_{min}(1,c,b) + \tfrac{1}{6} P_{max}(1,b,c) + \tfrac{1}{2} P_{max}(a+1,b,c), \nonumber\\
    &\tfrac{2}{3} - \tfrac{2}{3} P_{min}(1,c,a+b) + \tfrac{1}{3} P_{max}(1,a+b,c)   \Big). \nonumber
  \end{align}  
  Note that if a+b=n, then $P_{min}(a,b,c) := P_{max}(a,b,c) := 1$.
  
    \item $S(a,b,c)$ = continue \\
  $P_{min}(a,b,c) := \frac{1}{3} - \frac{1}{3} P_{max}(1,c,b) + \frac{1}{6} P_{min}(1,b,c) + \frac{1}{2} P_{min}(a+1,b,c)$
  
 $P_{max}(a,b,c) := \frac{1}{3} - \frac{1}{3} P_{min}(1,c,b) + \frac{1}{6} P_{max}(1,b,c) + \frac{1}{2} P_{max}(a+1,b,c)$\\

    %\medskip

    \item $S(a,b,c)$ = stop
    
  $P_{min}(a,b,c) := \frac{2}{3} - \frac{2}{3} P_{max}(1,c,a+b) + \frac{1}{3} P_{min}(1,a+b,c)$
  
 $P_{max}(a,b,c) := \frac{2}{3} - \frac{2}{3} P_{min}(1,c,a+b) + \frac{1}{3} P_{max}(1,a+b,c)$

     \medskip

\end{enumerate}

The iterations stop when $S(a,b,c)$ is either continue or stop for all positions.
With this method the optimal strategies for $n$ up to 20 were determined. Once an optimal strategy is found, the corresponding equations for the winning probabilities at each position are established and the system of equations is solved with a standard solver for linear equations.
The winning probabilities for the initial game situation (0, 0, 0) can be found in the two On-line Encyclopedia of Integer Sequences A387261 and A387262 \cite{A387261} for $n$ up to 15.

\section{Optimal strategy}

If player 1 needs only one more point, there is no strategy required because they toss only one more time. \\

In the case where the opponent needs only one more point, player 1 must go "all in", which means always continue. The proof is straightforward: The probability $P_{all}$ to win with the strategy  "all in" is
$$P_{all} = p_n \left (1 + \frac{1-p_n}{2} + (\frac{1-p_n}{2})^2 + \ldots \right ) = \frac{2p_n}{1+p_n}$$
with $p_n = 2^{-n}$ the probability to toss $n$ consecutive times heads.\\

If the player decides to bank an intermediate number of $x$ points, the probability $P_x$ to win is
$$P_x = \frac{2p_x}{1+p_x} \frac{p_{n-x}}{2} \left (1 + \frac{1-p_{n-x}}{2} + (\frac{1-p_{n-x}}{2})^2 + \ldots \right ) = \frac{2p_x}{1+p_x} \frac{p_{n-x}}{1+p_{n-x}}.$$

$P_x$ is smaller than $P_{all}$ for all $x$ and $n$: 
$$\frac{2p_x}{1+p_x} \frac{p_{n-x}}{1+p_{n-x}} < \frac{2p_n}{1+p_n} \Leftrightarrow$$
$$\frac{1 + p_x}{2p_x} \frac{1+p_{n-x}}{p_{n-x}} > \frac{1+p_n}{2p_n} \Leftrightarrow$$
$$(1+p_x)(1+p_{n-x})  > 1+p_n \Leftrightarrow p_x + p_{n-x} > 0.$$

Since breaking up a winning streak in two parts reduces the winning probability, the breaking up in more than two parts is even worse. Therefore "all in" is the best strategy in this case.\\

Remain the cases where both players need at least 2 more points.
A simple way to express the optimal strategy is to provide the number of coins a player should try to toss at any given game situation before converting the open points into saved points. Table~\ref{table:1} lists the optimal numbers for $n$ up to 20.\\

\begin{table}[h!]
\centering
\begin{tabular}[h]{|cc|c|c|c|c|c|c|c|c|c|c|c|c|c|c|c|c|c|c|c|}
\hline
& & \multicolumn{19}{c|}{Points Opponent needs} \\
& & 2 & 3 & 4 & 5 & 6 & 7 & 8 & 9 & {\rotatebox[origin=c]{90}{10}}& {\rotatebox[origin=c]{90}{11}}& {\rotatebox[origin=c]{90}{12}}& {\rotatebox[origin=c]{90}{13}}& {\rotatebox[origin=c]{90}{14}}& {\rotatebox[origin=c]{90}{15}}& {\rotatebox[origin=c]{90}{16}}& {\rotatebox[origin=c]{90}{17}}& {\rotatebox[origin=c]{90}{18}}& {\rotatebox[origin=c]{90}{19}}& {\rotatebox[origin=c]{90}{20}} \\
\hline
\parbox[t]{2mm}{\multirow{17}{*}{\rotatebox[origin=c]{90}{Points Player 1 needs}}}
& 2 & 2 & 2 & 1 & 1 & 1 & 1 & 1 & 1 & 1& 1& 1& 1& 1& 1& 1& 1& 1& 1& 1 \\
& 3 & 3 & 1 & 1 & 1 & 1 & 1 & 1 & 1 & 1& 1& 1& 1& 1& 1& 1& 1& 1& 1& 1 \\
& 4 & 2 & 2 & 2 & 1 & 1 & 1 & 1 & 1 & 1& 1& 1& 1& 1& 1& 1& 1& 1& 1& 1 \\
& 5 & 2 & 2 & 2 & 2 & 1 & 1 & 1 & 1 & 1& 1& 1& 1& 1& 1& 1& 1& 1& 1& 1 \\
& 6 & 2 & 2 & 2 & 2 & 2 & 1 & 1 & 1 & 1& 1& 1& 1& 1& 1& 1& 1& 1& 1& 1 \\
& 7 & 3 & 2 & 2 & 2 & 2 & 2 & 1 & 1 & 1& 1& 1& 1& 1& 1& 1& 1& 1& 1& 1 \\
& 8 & 3 & 2 & 2 & 2 & 2 & 2 & 2 & 1 & 1& 1& 1& 1& 1& 1& 1& 1& 1& 1& 1 \\
& 9 & 3 & 3 & 2 & 2 & 2 & 2 & 2 & 2 & 1& 1& 1& 1& 1& 1& 1& 1& 1& 1& 1 \\
& 10 & 3 & 2 & 2 & 2 & 2 & 2 & 2 & 2 & 2& 1& 1& 1& 1& 1& 1& 1& 1& 1& 1 \\
& 11 & 3 & 3 & 2 & 2 & 2 & 2 & 2 & 2 & 2 & 2 & 1 & 1 & 1 & 1 & 1 & 1 & 1 & 1 & 1\\
& 12 & 3 & 3 & 2 & 2 & 2 & 2 & 2 & 2 & 2 & 2 & 2 & 1 & 1 & 1 & 1 & 1 & 1 & 1 & 1\\ 
& 13 & 4 & 3 & 3 & 2 & 2 & 2 & 2 & 2 & 2 & 2 & 2 & 2 & 1 & 1 & 1 & 1 & 1 & 1 & 1\\ 
& 14 & 4 & 3 & 2 & 2 & 2 & 2 & 2 & 2 & 2 & 2 & 2 & 2 & 2 & 1 & 1 & 1 & 1 & 1 & 1\\ 
& 15 & 4 & 3 & 3 & 2 & 2 & 2 & 2 & 2 & 2 & 2 & 2 & 2 & 2 & 2 & 1 & 1 & 1 & 1 & 1\\ 
& 16 & 4 & 3 & 3 & 2 & 2 & 2 & 2 & 2 & 2 & 2 & 2 & 2 & 2 & 2 & 2 & 1 & 1 & 1 & 1\\ 
& 17 & 4 & 3 & 3 & 3 & 2 & 2 & 2 & 2 & 2 & 2 & 2 & 2 & 2 & 2 & 2 & 2 & 1 & 1 & 1\\ 
& 18 & 4 & 3 & 3 & 3 & 2 & 2 & 2 & 2 & 2 & 2 & 2 & 2 & 2 & 2 & 2 & 2 & 2 & 1 & 1\\ 
& 19 & 4 & 3 & 3 & 3 & 2 & 2 & 2 & 2 & 2 & 2 & 2 & 2 & 2 & 2 & 2 & 2 & 2 & 2 & 1\\ 
& 20 & 4 & 3 & 3 & 3 & 2 & 2 & 2 & 2 & 2 & 2 & 2 & 2 & 2 & 2 & 2 & 2 & 2 & 2 & 2\\
\hline 
\end{tabular}
\caption{Number of coins a player should toss to optimize their winning chances.}
\label{table:1}
\end{table}

For example if both players need 4 more points to win, the player whose turn it is, should continue if they toss heads in the first round, but they should stop if they toss heads a second time, thus converting two open points.\\

\section*{Acknowledgements}

The author would like to thank Jinyuan Wang for contributing the iterative method to determine the optimal strategy; Lucas Wey Hacker for drawing my attention to the paper by Glenn et al.; Kester Habermann for coding support; and Tatiana Tatarenko, Misha Lavrow and Denizalp Goktas for their valuable advice.

\bibliography{references} 
\bibliographystyle{plain} %ieeetr

\end{document}